\documentclass[letterpaper, 12pt]{article}
   
\usepackage{amsmath, amssymb,amsthm,tikz}
\usepackage{enumerate}
\usepackage[title]{appendix}
\usepackage{fullpage}

\usepackage{color}
\usepackage{colortbl}
\definecolor{mygray}{gray}{1}
\definecolor{mygray2}{gray}{.9}
\definecolor{mygray3}{gray}{.8}

\usepackage[bookmarks]{hyperref}

\def\half{{\frac{1}{ 2}}}
\def\reci#1 {{\frac{1}{#1}}}

\gdef\Rk{{ R_k}}       

\def\proofend
{\ifhmode
\unskip\nobreak\hfill\vrule height5pt width4pt depth2pt\medskip\fi
\ifmmode\eqno{\vrule height5pt width4pt depth2pt}\fi}

\def\Rd{}

\gdef\Ta{{t}}
\gdef\da{{d}}

\gdef\sa{{s}}
\gdef\pa{{s}}

\def\beq#1{\begin{equation}\label{#1}}
\def\eeq{\end{equation}}

\newtheorem{theorem}{Theorem}[section]

\newtheorem{conjecture}[theorem]{Conjecture}

\newtheorem{corollary}[theorem]{Corollary}

\newtheorem{lemma}[theorem]{Lemma}

\newtheorem{claim}[theorem]{Claim}

\theoremstyle{definition}
\newtheorem{definition}[theorem]{Definition}

\newtheorem{construction}[theorem]{Construction}
\newtheorem{question}[theorem]{Question}

\def\Question#1{\medskip
\begin{question}
#1
\end{question}
}

\def\Proof#1. {\bigskip\noindent{\Rd\sc Proof #1.~}}

\gdef\proof{\bigskip\noindent{\Rd\sc Proof.~}}

\def\text#1{\qquad\mbox{#1}\qquad}

\def\Text#1{\quad\mbox{#1}\quad}

\def\eps{\varepsilon}
\def\ph{\varphi}

 \def\Gn{{ G_n}}

\def\Sn{{ S_n}}

\def\Tnr{{\Rd T_{n,r}}}

\def\RT{{\rho\tau}}
\gdef\RTP{{\bf RT}}
\gdef\RTT{{\bf RT_d}}
\def\RTS{{\overline{\rho\tau}}}
\def\RTI{{\underline{\rho\tau}}}

\def\R{{\bf R}}
\def\Q{{\bf Q}}

\def\remark{{\medskip\noindent \textbf{Remark:}} }
\def\Remark{{\noindent \textbf{Remark:}} }

\def\w#1{\omega(#1)}
\def\q#1{\Q\left(#1, f\right)}
\def\qn#1{\Q\left(#1, f(n)\right)}
\def\oq#1{\Q\left(#1, \frac{n}{\w{n}}\right)}
\def\oqd#1{\Q\left(#1, \dfrac{n}{\w{n}}\right)}
\def\nq{\dfrac{n}{q}}
\def\nq{n}

\gdef\td#1{\frac{1}{2}\left(1-\frac{1}{#1}\right)}
\gdef\on{ o\left(n^2\right)}

\let\originalleft\left
\let\originalright\right
\renewcommand{\left}{\mathopen{}\mathclose\bgroup\originalleft}
\renewcommand{\right}{\aftergroup\egroup\originalright}

\title{Phase transitions in  Ramsey-Tur\'an theory}

\author{J\'ozsef Balogh \footnote{
    Department of Mathematics, University of Illinois, Urbana, IL 61801, USA, {\tt jobal@math.uiuc.edu}.
    Research is partially supported by NSF CAREER Grant DMS-0745185 and Arnold O. Beckman Research Award (UIUC Campus Research Board 13039) and Marie Curie FP7-PEOPLE-2012-IIF 327763.}
\and Ping Hu \footnote{Department of Mathematics, University of Illinois, Urbana, IL 61801, USA, {\tt pinghu1@math.uiuc.edu.}}
\and Mikl\'os Simonovits \footnote{R\'enyi Institute, Budapest, Hungary, {\tt miki@renyi.hu}.
 Research is partially supported by
OTKA grant 101536.}
}

\date{}

\begin{document}
\maketitle

\begin{abstract}
 Let $f(n)$ be a function and $H$ be a graph.  Denote by
$\RTP(n,H,f(n))$ the maximum number of edges of an $H$-free graph on
$n$ vertices with independence number less than $f(n)$.  Erd\H os and
S\'os~\cite{ES} asked if $\RTP\left(n, K_5, c\sqrt{n}\right) = \on$
for some constant $c$.  We answer this question by proving
the stronger $\RTP\left(n, K_5, o\left(\sqrt{n\log n}\right)\right) = \on$. It is known that
$\RTP \left(n, K_5, c \sqrt{n\log n} \right) = n^2/4+\on $ for $c>1$,
so one can say that $K_5$ has a Ramsey-Tur\'an
{\em phase transition} at $c\sqrt{n\log n}$.  We extend this result
to several other $K_\pa$'s and functions $f(n)$,
determining many more phase transitions.  We 
formulate several open problems, in particular, whether variants of
the Bollob\'as-Erd\H os graph exist to give good lower bounds on
$\RTP\left(n, K_\pa, f(n)\right)$ for various pairs of $\pa$ and
$f(n)$. Among others, we use Szemer\'edi's Regularity Lemma and the
Hypergraph Dependent Random Choice Lemma. We also present a short
proof of the fact that $K_\pa$-free graphs with small independence
number are sparse: on $n$ vertices have $o(n^2)$ edges.

\medskip
Keywords: Ramsey, Tur\'an, independence number, dependent random choice
 \end{abstract}


\section{Introduction}
{
{\bf Notation.} In this paper we shall consider only simple graphs, i.e., graphs without loops
and multiple edges. As usual, $\Gn$ will always denote a graph on
$n$ vertices. More generally, in case of graphs the (first) subscript
will always denote the number of vertices, for example $K_\pa$ is the
complete graph on $\pa$ vertices, and $\Tnr$ is the $r$-partite Tur\'an
graph on $n$ vertices, i.e., the complete $r$-partite graph on $n$
vertices with class sizes as equal as possible. Given a graph $G$, we
use $e(G)$ to denote its number of edges, and use
$\alpha(G)$ to denote its independence number. Given a subset $U$ of the vertex set of  $G$, we use $G[U]$ to denote the subgraph of $G$ induced by $U$.

  In this paper all logarithms are base $2$; $\w{n}$ denotes an
  arbitrary function tending to infinity 
slowly enough so that all calculations we use go through. 
 Whenever we write
  that  ``$\omega(n)\to\infty$ slowly", we mean that
  the reader may choose an arbitrary $\omega(n)\to\infty$, the assertion
  will hold, and the more slowly $\omega(n)\to\infty$ the stronger the
  assertion, i.e., the theorem is. In the proofs, we shall assume that $\w{n}=o(\log\log\log n)$.
  In our cases, if we prove some theorems for such functions $\w{n}$, then these results remain valid for larger functions as well.
To simplify the formulas, we shall often omit the floor and
  ceiling signs, when they are not crucial.

S\'os~\cite{SosCalgary} and Erd\H os and S\'os~\cite{ES} defined the
following `Ramsey-Tur\'an' function:

\begin{definition}
Denote by $\RTP(n,H,f(n))$ the maximum number of edges of an $H$-free
graph on $n$ vertices with independence number less than $f(n)$.  
\end{definition}

%
%
Sometimes we want to study the case 
when
the  bound on the independence number $f(n)$ is $o(g(n))$.  Formally $o(g(n))$ is not a function, we shall consider
$\RTP\big(n,H, o\big(g(n)\big)\big)$ as $\RTP\big(n,H, g(n)/\w{n}\big)$ where $\w{n}$ is an arbitrary function tending to
infinity  (slowly).

More formally,  if $\RTP(n,H, f(n)) \le cn^2+o(n^2)$ for every $f(n) = o(g(n))$, then  we write $\RTP\big(n, H, o\big(g(n)\big)\big) \le cn^2+o(n^2)$. If $\RTP(n, H, f(n)) \ge cn^2+\on$ for some $f(n) = o(g(n))$,
then we write $\RTP\big(n, H, o\big(g(n)\big)\big) \ge cn^2+\on$. If we can show both  $\RTP\big(n, H, o\big(g(n)\big)\big)  \le cn^2+o(n^2)$ and $\RTP\big(n, H, o\big(g(n)\big)\big) \ge cn^2+\on$, then we write $\RTP\big(n, H, o\big(g(n)\big)\big)  = cn^2+\on$.


 We are interested in the asymptotic
  behavior of $\RTP(n,H,f(n))$, i.e., if $\RTP(n,H, f(n)) = cn^2+o(n^2)$, then we are more interested in the constant $c$ than in the behavior of $\on$.
  
  \begin{definition}
Let \begin{align*}
\RTS(H, f) = \limsup_{n\to\infty} \frac{\RTP\big(n, H, f(n)\big)}{{n^2}} \Text{and}
\RTI(H, f) = \liminf_{n\to\infty}  \frac{\RTP\big(n, H, f(n)\big)}{{n^2}}.
\end{align*}
If $\RTS(H, f) = \RTI(H, f)$, then we write $\RT(H, f) = \RTS(H, f) =
\RTI(H, f)$, and call $\RT$ the {\em Ramsey-Tur\'an density} of $H$
with respect to $f$, $\RTS$ the upper, $\RTI$ the lower {\em
  Ramsey-Tur\'an densities}, respectively.
\end{definition}
\Remark It is easy to see that $\RT(H, f) = c$ is equivalent to
$\RTP(n,H,f(n)) =cn^2+\on$. {Here we define $\RTS$ and $\RTI$ simply because for an arbitrary  function $f$, the limit $\RT$ may not exist. For reasonable functions $f$ like most functions considered in this paper, we can show $\RT$ exists.}
When we write $\RT(H, f)$, we use $f$ instead of $f(n)$, since
$\RT(H,f(n))$
 would suggest that this constant depends on $n$. If
however, we write something like $\RT\left(H,c\sqrt{n\log n}\right)$, that is
(only) an abbreviation of $\RT(H,f)$, where $f(n)=c\sqrt{n\log n}$,
(see e.g. Theorem~\ref{sharp_k5}).  So, even when we write $\RT(H,
f(n))$, we are treating $f(n)$ as a function, which means $\RT(H,
f(n))$ does not depend on $n$.
\medskip

 We try to understand that given a graph $H$ and  a very large $n$, when do we observe crucial drops in the
  value of $\RTP(n,H,m)$ while $m$ is changing (continuously) from $n$ to $2$? In other words, we try to
understand
when and how the asymptotic behavior of $\RTP(n,H,f(n))$ changes sharply
  when we replace $f$ by a slightly smaller $g$.  

\begin{definition}[Phase Transition]\label{PhTDef}

Given a graph $H$ and two functions $g(n) \le f(n)$, we shall
say that $H$ has a {\em phase transition} from $f$ to $g$ if
$ \RTS(H,g) < \RTI(H,f)$.

Given a function $\ph(n)\to 0$,
we shall say that $H$ has a {\em
  $\ph$-phase-transition} at $f$ if $H$ has a phase transition from
$f$ to $\ph f$. If $H$ has a $\ph$-phase-transition at $f$ for every
$\ph$ tending to $0$, then we shall say that $G$ has a {\em strong
  phase transition} at $f$. Let $\ph_{\eps}(n) =
2^{-\log^{1-\eps} n}$.
If there exists an $\eps>0$ for which $H$
has a $\ph_{\eps}$-phase-transition at $f$,
 then we shall say that $H$ has a {\em
  weak phase transition} at $f$.\footnote{The strange function $2^{\log^{1-\eps}n}$ is somewhere
  ``halfway'' between $\log n$ and $n^c$.  }
\end{definition}
The Ramsey-Tur\'an  theory  is very complex,
with many open questions. 
Here we focus on the case when $H$ is a clique.
Erd\H{o}s and S\'os~\cite{ES} determined $\RTP(n,K_{2r+1},o(n))$.

\begin{theorem}
\label{ES:RT}
For every positive integer $r$, 
$$\RTP\big(n,K_{2r+1},o(n)\big)=\td r n^2+o(n^2).$$
\end{theorem}

The meaning of Theorem~\ref{ES:RT} is that the Ramsey-Tur\'an density
of $K_{2r+1}$ in this case is essentially the same as the Tur\'an
density $\half(1-1/r)$ of $K_{r+1}$. It also shows that $K_5$ has a strong phase transition at
$n$, since, by Tur\'an's Theorem, we have $\RT( K_5, n) =3/8$.  
 In~\cite{ES}, Erd\H{o}s and S\'os
proved that $\RTP\left(n, K_5, c\sqrt{n}\right) \le  n^2/8 + \on$ for every
$c>0$.\footnote{
It follows from the proof of Erd\H{o}s and S\'os and the result of Ajtai, Koml\'os and Szemer\'edi~\cite{AKS} on the Ramsey number $\R(3,n)$ that $\RTP\left(n, K_5, o\left(\sqrt{n\log n}\right)\right) \le  n^2/8 + \on$.}
 They also asked if $\RTP\left(n, K_5, c\sqrt{n}\right) = \on$ for some
$c>0$.  One of our main results, Theorem~\ref{sharp_k5}, states { $\RTP(n, K_5,$ $ o\left(\sqrt{n\log n}\right)) = o(n^2)$}, answering this question. This result together with \eqref{eq:k5} shows that $K_5$ has a strong
phase transition at $c\sqrt{n\log n}$ for every $c>1$.
Actually, every $K_\pa$ with $\pa>2$ has a strong phase
  transition at $n$.
More generally, given a graph $G$, 
if $\chi(G)>2$ and $G$ has an edge $e$ such that $\chi(G-e)<\chi(G)$, then $G$ has a strong phase transition at $n$.
On the other hand, let $K_\pa(a_1,\dots, a_\pa)$ be the complete $\pa$-partite graph with class sizes $a_1,\dots, a_\pa$.
Simonovits and S\'os~\cite{SimSosRT} showed that if $\pa<a\le b$, then $$\RT(K_{\pa+1}(a,b,\dots,b), o(n)) = \RT(K_{\pa+1}(a,b,\dots,b), n) = \td{\pa},$$ 
which means that $K_{\pa+1}(a,b,\dots,b)$ does not have a strong phase
  transition at $f(n)=n$.

Clearly, a  $\phi$-phase-transition implies a $\phi'$-phase-transition if $\phi$ tends to $0$ more slowly than $\phi'$, and in particular
a strong phase transition implies a weak phase transition.

The rest of the paper is organized as follows. In
Section~\ref{sec:history} we provide additional history of Ramsey-Tur\'an
type problems.  Our aim in general is to determine the phase
transitions for cliques, we state 
our new results in Section~\ref{sec:results}.  In
Section~\ref{sec:tools}
 we provide the main tools for our proofs: the
Dependent Random Choice Lemma and the Hypergraph Dependent
Random Choice Lemma.  We prove our main results in
Sections~\ref{sec:proof_drc} and \ref{sec:proof_hydrc}, 
 and
Section~\ref{sec:last} contains some 
open problems.



\section{History}\label{sec:history}
\subsection{Ramsey Numbers}
In order to better understand the Ramsey-Tur\'an theory, we need some results from Ramsey theory.
There are several constructions providing lower
  bounds on Ramsey-Tur\'an functions  that are based on constructions
 corresponding to some ``simple, small Ramsey numbers''. 
   Let $\R(\Ta,m)$ be the
Ramsey number: the minimum $n$ such that every graph $\Gn$ on $n$
vertices contains a clique $K_\Ta$ or an independent set of size $m$.

Unfortunately, we do not know Ramsey functions very well.  The case
$\Ta=3$ is well-understood.
The bound $\Q(3,n) =\Theta\left( \sqrt{n\log n} \right)$ was proved by
Ajtai, Koml\'os, Szemer\'edi~\cite{AKS} and Kim~\cite{kim}. The best known
quantitative estimates were proved by Shearer~\cite{Shearer}, Pontiveros, Griffiths, Morris~\cite{PGM} and Bohman, Keevash~\cite{BK2}. The bounds are
  \begin{align}
\left(1\middle/4 - o(1)\right)\frac{ m^2 }{ \log m} \le \R(3,m)  \le \left(1 + o(1)\right) \frac{m^2} { \log m}. \label{r3t} 
\end{align}


For $\Ta\ge 4$ we have only
\beq{eq:ramsey}\Omega\left(\frac{m^{(\Ta+1)/2}}{(\log m)^{(\Ta+1)/2-1/(\Ta-2)}}\right)\le   \R(t,m)\le O\left(\frac{m^{\Ta-1}}{(\log m)^{\Ta-2}}\right),\eeq
where the upper bound follows from Ajtai, Koml\'os and Szemer\'edi~\cite{AKS} and the
lower bound follows from Bohman and Keevash~\cite{BK}. It is conjectured that the upper bound is sharp up to some $\log m$-power factors. 

We define the `inverse' function $\Q(\Ta,n)$ of $\R(\Ta,m)$, i.e., the
minimum independence number of $K_t$-free graphs on $n$ vertices.
It is an inverse function
 in the sense that if 
$\R(\Ta, m) = n$, then $\Q(\Ta,n) = m$.
For example, $\Q(2, n) = n, \Q(3, n) = \Theta(\sqrt{n\log n})$ and 
$\Omega (n^{1/3} \log^{2/3} n)   \le \Q(4, n)=O(n^{2/5} \log^{4/5} n).$
 In general, for $t\ge 3$, we know from \eqref{r3t} and \eqref{eq:ramsey} that
\begin{align}
   \left(\frac{1}{\sqrt{2}}-o(1)\right)\sqrt{n\log n} \le & \Q(3,n)  \le  \left(\sqrt{2} + o(1)\right) \sqrt{ n \log n }\label{beta}\\
   \Omega\left(   n^{\frac{1}{\Ta-1}} (\log n)^{\frac{\Ta-2}{\Ta-1}}   \right)  \le& \Q(\Ta,n)\le  O\left( n^{\frac{2}{t + 1}} (\log n)^{ 1- \frac{2}{(\Ta-2)(\Ta+1)}} \right)\label{Q}.
   \end{align}

There are some
famous conjectures on $\R(\ell, n)$. We state three of them, with increasing strength.

\begin{conjecture}\label{conj:Qsn}
(a)  (Folklore) For every integer $\ell\ge 3$, $\R(\ell-1, n) =o(\R(\ell, n))$ as $n\to\infty$.\\
(b) For every integer $\ell \ge 3$, there exist $\vartheta = \vartheta(\ell) >0 $ and $N= N(\ell) > 0$ such that if $n>N$, then
\beq{eq:Qsn}
\R(\ell -1, n) \le \frac{\R(\ell, n)}{n^{\vartheta}}. 
\eeq 
(c)  For some constant $\gamma=\gamma(\Ta)$,
$$\Q(\Ta,n)\approx \root{ \Ta-1} \of n\log^{\gamma}n,$$
or at least
$$ \root{ \Ta-1} \of n<\Q(\Ta,n)< \root{ \Ta-1} \of n\log^{\gamma}n.$$
\end{conjecture}

We know { from \eqref{eq:ramsey}} 
that Conjecture~\ref{conj:Qsn} is true for $\ell=3,4$, but 
 for larger $\ell$'s we are 
very far from proving what is conjectured.

\subsection{History  of Ramsey-Tur\'an Theory}

 Let $H_{k,\ell}$ denote a ``Ramsey'' graph on $k$ vertices not
 containing $K_{\ell}$, having  the minimum possible independence number under this condition.
 The graph $H_{k,\ell}$ is sparse,
   i.e., it has $o(k^2)$ edges, see Theorem~\ref{thm:ramsey_sparse}.
 For Theorem~\ref{ES:RT}, Erd\H{o}s and S\'os~\cite{ES} used $H_{n/r,
   3}$ to construct a  graph $\Sn$ to provide the lower bound on 
 $\RT(K_{2r+1},
 o(n))$.  Their idea was that when a Ramsey-graph
 $H_{n/r,3}$ is placed into each class of a Tur\'an graph $\Tnr$, we
 get a $K_{2r+1}$-free graph sequence $\{\Sn\}$ with \beq{LowerRTGen}
 e(\Sn)\approx
 e(\Tnr)\text{and}\alpha(\Sn)=\alpha(H_{n/r,3})=o(n).  \eeq It is 
trivial to generalize this idea to give a
 lower bound on $\RT(K_{rs+1},o(n))$.

\begin{construction}[Extended/Modified Erd\H os-S\'os Construction]\label{LowerRT}
~\\Let $k=\left\lfloor n/r \right\rfloor$,
take a Tur\'an graph $\Tnr$ with $r$ classes 
and place an $H_{k, t+1}$ into each of its classes.
\end{construction}

It is easy to see that this graph is
$K_{rt+1}$-free, hence 
\beq{lowerBound} \RTP(n, K_{rt+1}, \alpha(H_{n/r, 
  t+1}) ) \ge \td r n^2+\on.  \eeq If Conjecture~\ref{conj:Qsn}~(b) is true for
$\ell = t+1$, then $\RT(K_{rt +1}, \alpha(H_{n/r,  t+1}) )$ exists and
\eqref{lowerBound} is sharp, see
Theorem~\ref{thm:Kpq}.
\medskip

Szemer\'edi~\cite{szemrt}, using  an earlier, weaker form of his regularity lemma~\cite{RegLemma}, proved
$\RTS(K_4,$ $ o(n)) \le 1/8$. Bollob\'as and Erd\H os~\cite{BE}
constructed the so-called Bollob\'as-Erd\H os graph, one of the most
important constructions in this area, 
 that  shows that $\RTI(K_4, o(n)) \ge
1/8$.  
Indeed, the Bollob\'as-Erd\H os graph on $n$ vertices is $K_4$-free,
with $(\reci 8 + o(1))n^2$ edges and independence number $o(n)$.  Later,
Erd\H os, Hajnal, S\'os and Szemer\'edi~\cite{EHSSz}  extended
  these results, determining  $\RTP(n, K_{2r},o(n))$:


\begin{theorem}\label{K2r:RT}
$$\RTP(n,K_{2r},o(n))=\frac{3r-5}{6r-4}n^2+\on.$$
\end{theorem}

The lower bound is provided by their 
generalization of
the Bollob\'as-Erd\H os graph:

\begin{construction}\label{LowerBE}

Fix $h = \left\lfloor \frac{4n}{3r-2}\right\rfloor$ and
$k=\left\lfloor \frac{3n}{3r-2}\right\rfloor$. Let $B_h$ be
a Bollob\'as-Erd\H os graph on $h$ vertices.  We take
 vertex-disjoint copies of  a
$B_h$ and a Tur\'an graph $T_{n-h, r-2}$,  and join each vertex of $B_h$ to
each vertex of $T_{n-h,r-2}$, and place an $H_{k,3}$ into each class
of $T_{n-h,r-2}$. 
\end{construction}

Here $h$ was chosen to maximize the number of
  edges, which is equivalent with making the degrees (almost) equal.
It is easy to see that this graph is $K_{2r}$-free. Since $\alpha
(B_h)=o(n)$ and $\alpha (H_{k,3}) = o(n)$, it gives the lower bound of
Theorem~\ref{K2r:RT}. \medskip

Recall that $\Q(\Ta ,n)$ is the minimum independence number of
$K_\Ta$-free graphs on $n$ vertices.  So we have $\Q(\Ta +1, n/r) =
\alpha(H_{n/r, \Ta+1})$, and we can write \eqref{lowerBound} as
\beq{lowerQ} \RTI\left(K_{r\Ta+1}, \Q\left( \Ta+1,
\frac{n}{r}\right)\right ) \ge\td{r}= \frac{r-1}{2r}.  \eeq In
particular, we get the following sharpening of the lower bound of
Theorem~\ref{ES:RT}:
\begin{equation}\label{k5upperb}
\RTI\left(K_{2r+1}, \Q\left(3, \frac{n}{r}\right)\right)\ge 
\td r =
\frac{r-1}{2r}. 
\end{equation}

Combining \eqref{k5upperb} and \eqref{beta},
we have the following relation.\footnote{Essentially this appears in Erd\H os-S\'os~\cite{ES}.}
For any $c>1,$ 
\begin{equation}\label{k5upperC}
\RTI\left(K_6,c\sqrt{n\log n}\right)\ge \RTI\left(K_5,c\sqrt{n\log n}\right)\ge \frac{1}{4}.
\end{equation}

Erd\H{o}s, Hajnal,
Simonovits, S\'os and Szemer\'edi~\cite{EHSSSzB} considered $K_6$ and proved the following theorem. 

\begin{theorem}\label{thm:SS}

If $\w{n}\to\infty$, then
\begin{align*}
 \RTP\left(n, K_{6}, \Q\left(3,\frac{n}{\w{n}}\right)\right) \le \frac{n^2}{6}+\on.
\end{align*}
\end{theorem}

In the last years, many important, new results were proved 
 on $\RT(K_4,o(n))$. 
      Sudakov proved that
$\RT\left(K_4, n2^{-\w{n}\sqrt{\log n}}\right) = 0$, which  is a special case of his  more general theorem~\cite{sud_rt}:

\begin{theorem}\label{thm:sudk2t}
Let $t\ge 2$ and $\w{n}\to\infty$. If $g(n) = \Q\left(t, n2^{-\w{n}\sqrt{\log n}}\right)$, 
  then $\RT(K_{2t}, g)= 0$.
\end{theorem}

Recently, by finding good quantitative
estimates for the relevant parameters of the  Bollob\'as-Erd\H os
graph, Fox, Loh and Yufei Zhao~\cite{FLZ} proved that $\RTI\left(K_4, n2^{-o\left(\sqrt{\frac{\log
      n}{\log \log n}}\right)}
\right)\ge 1/8$, complementing Sudakov's
result on $K_4$.

Ramsey-Tur\'an problems with independence
number $\Q(\Ta, f(n))$ were also studied earlier in
a somewhat different way. Given an integer 
$d\ge 2$, define the {\em $d$-independence number} $\alpha_d(G)$ of $G$
to be the maximum size of a vertex set $S$ for which $G[S]$ contains no $K_d$.
For example, the independence number $\alpha(G)$ of $G$ is $\alpha_2(G)$.
Denote by $\RTT(n,H,f(n))$ the maximum number of edges of an $H$-free graph on
$n$ vertices with $d$-independence number less than $f(n)$.  It is easy to see
that $\alpha(\Gn)<\Q(\da, f(n))$ implies $\alpha_\da(\Gn) < f(n)$, so
$\RTP(n,H,\Q(\da, f(n)))\le\RTT(n,H, f(n))$. Therefore an upper bound on 
  $\RTT(n,H, f(n))$ is also an upper bound on $\RTP(n,H,\Q(\da,f(n)))$.
Erd\H{o}s, Hajnal, Simonovits, S\'os and Szemer\'edi~\cite{EHSSSzB} gave an
upper bound on $\RTT(n,K_\pa, o(n))$, that  implies the following theorem.

\begin{theorem}\label{erdsosind}

For any function $\w{n}$ tending to infinity, if $2\le \Ta<\sa$, then
\beq{eq:ind}
\RTP\left(n,K_\sa,\Q\left(\Ta,\frac{n}{\w{n}}\right)\right)\le \frac{\sa-\Ta-1}{2\sa-2}n^2+\on.
\eeq
\end{theorem}

 Lower bounds on $\RTT(n, K_\pa, o(n))$
 were provided by constructions of Balogh and Lenz~\cite{BL1, BL2}.
  Unfortunately, a lower bound on $\RTT(n,H,f(n))$ provides no
 lower bound on $\RTP(n,H,$ $\Q(d,f(n)))$. For example, Balogh and Lenz
 \cite{BL2} gave a construction showing that $\RTP_3(n, K_5,$ $ f(n) )
 \ge  n^2/16+\on$ for some $f(n)=o(n)$; on the other hand,
 Theorem~\ref{sharp_k5} implies that 
 $\RTP(n,K_5,\Q(3,f(n))) = \on$ for any $f(n)=o(n)$.

\section{New Results}\label{sec:results}
First,
we show that
$K_{\ell}$-free graphs with small  independence number are sparse.

\begin{theorem}\label{thm:ramsey_sparse}

Let $\ell\ge3$ be an integer and $s = \left\lceil \ell/2\right\rceil$. 
 Fix a positive constant  $c < \frac{1}{s(s-1)}$.
Let $G_{n,\ell}$ be a graph on $n$ vertices not containing $K_{\ell}$. 
$$\text{If}
\alpha(G_{n,\ell}) < \Q(\ell,n)n^c, \text{then} e(G_{n, \ell})=o(n^2).$$
\end{theorem}

%
%
%
%
%
%
%
One of our main results is the following.
\begin{theorem}\label{sharp_k5}
If $\w{n}\to \infty$, then
\begin{equation}\label{ccc}
\RTP\left(n, K_5, \frac{\sqrt{n\log n}}{\w{n}}\right)\le \frac{n^2}{\sqrt[4]{\w{n}}}=\on.
\end{equation}
\end{theorem}

Here, by Construction~\ref{LowerRT}, $\sqrt{n\log n}/\w{n}$
 is sharp in the sense that 
 \beq{eq:k5}
\quad\RTI\left(K_5,c\sqrt{n\log n}\right)\ge \td 2=\frac{1}{4} \text{for any} c>1. \eeq

We generalize Theorem~\ref{sharp_k5} from  $K_5$ to many other $K_s$.

\begin{theorem}\label{thm:pq-1}
 Suppose  $p\ge 3$ and $q\ge 2$. Let $\w{n}\to \infty$ and $f(n) = n2^{-\w{n}\log^{\frac{q-2}{q-1}} n}.$ \\
 (a) If there exists a constant $\vartheta>0$ such that for every $n$ sufficiently large,  we have \newline
$\R\big(p-1$, $\qn{p}\big)< n^{1-\vartheta}$, then  
 \begin{equation}\label{eq:pq-1} 
   \RTS\left(K_{pq-1}, \q{p}\right) \le  \frac{1}{2}\left(1- \frac{1}{q-1}\right).
  \end{equation}
(b) For $1\le i \le p-1$, let 
$t =\left\lfloor\frac{pq-i-1}{q-1}\right\rfloor \ge p$.
If $$\Q\left(t+1, \frac{n}{q-1}\right) \le g(n) \le \Q\left(p, f(n)\right),
\Text{then}
\RT\left(K_{pq-i}, g\right) =  \frac{1}{2}\left(1- \frac{1}{q-1}\right).$$
\end{theorem}

 We extend Theorem~\ref{thm:sudk2t} from $K_{2t}$ to $K_{pq}$, where $q$ replaces $2$ and $t$ replaces $p$. 
 Theorem~\ref{thm:Kpq}~(a) can be compared to
  Theorem~\ref{thm:pq-1}~(a), where similar statement was proved for
 $K_{pq-1}$ and a slightly larger $f(n)$.

\begin{theorem}\label{thm:Kpq}
   Suppose $p\ge2$ and $q\ge 2$. Let $\w{n}\to\infty$ and 
$f(n) = n2^{-\w{n}\log^{1-1/q }n}$, then\\
(a) \begin{equation}\label{eq:Kpq}
 \RTS\left(K_{pq}, \Q\left(p,f\right)\right) \le \frac{1}{2}\left(1- \frac{1}{q-1}\right).
 \end{equation}
(b) {For $0\le i\le p-1$, let $t=\left\lfloor\frac{pq-i-1}{q-1}\right\rfloor \ge p$. 
If $$\Q\left(t+1, \frac{n}{q-1}\right) \le g(n) \le \Q\left(p, f(n)\right),
\Text{then}
\RT\left(K_{pq-i}, g\right) =  \frac{1}{2}\left(1- \frac{1}{q-1}\right).$$}
\end{theorem}

{We generalize Theorem~\ref{thm:SS} from $K_6$ to all even cliques.}
\begin{theorem}\label{thm:K2t}

Let $\w{n}\to\infty$ and $f(n) = n2^{-\w{n}\log^{1-1/q }n}$. 
If $2t\le pq$ 
and $\oq{t} \le  \Q\left(p,f(n)\right) $, then
  \[ 
  \RTS\left(K_{2t}, \Q\left(t,\frac{n}{\w{n}}\right)\right) \le 
 \frac{1}{2}\left(1- \frac{1}{t}\right)     \left(1- \frac{1}{q-1}\right).
  \]

\end{theorem}

\bigskip

Using Theorem~\ref{erdsosind} and
\eqref{lowerQ}, we get the following corollary.

\begin{corollary}\label{Cor:erdsosind} 

Suppose $p, q\ge 2$. Let $\w{n}\to\infty$.
 If $\Q\left(p+1,
n/q \right)\le \Q\left(p,\frac{n}{\w{n}}\right)$,
then $$\RT\left(K_{pq+1}, \Q\left(p,\frac{n}{\w{n}}\right)\right) =   \frac{1}{2}\left(1- \frac{1}{q}\right).$$
 \end{corollary}

\proof
The upper bound follows from
Theorem~\ref{erdsosind} with $s=pq+1$ and $t=p$. 
The lower bound follows from \eqref{lowerQ} with $r=q$ and $\Ta=p$:
$$\RTI\left(K_{pq+1}, \Q\left(p+1, \frac{n}{q}\right) \right)\ge    \frac{1}{2}\left(1- \frac{1}{q}\right). \proofend$$ 

Now we are ready to
find phase transitions.  

Note that $\RTS\left(K_\pa,\Q\left(\Ta,n\right)/\w{n}\right) \le \RTS\left(K_\pa,\Q\left(\Ta,n/\w{n}\right)\right)$, so Theorem~\ref{erdsosind} gives an upper bound on $\RTS\left(K_\pa,\Q\left(\Ta,n\right)/\w{n}\right)$.
 Construction~\ref{LowerRT} provides $K_\pa$-free graphs with many edges and small independence number, 
 giving a lower bound on $\RTI\left(K_\pa,\Q\left(\Ta,n\right)\right)$.
Using these two results,  we give
conditions  on $\pa$ and 
 $\Ta$ under which 
$\RTS\left(K_\pa,\Q\left(\Ta,n\right)/\w{n}\right) < \RTI\left(K_\pa,\Q\left(\Ta,n\right)\right)$ for any $\w{n}$ tending to infinity, i.e.,
$K_\pa$ has a
strong phase transition at $\Q(\Ta,n)$.  

\begin{theorem}\label{thm:pt}

 If $\pa-1=
r(\Ta-1)+\ell$ with $0\le \ell < \Ta-1, \ell < r$ and $2\le \Ta<\pa$, then
$K_\pa$ has a strong phase transition at $f(n)=\Q(\Ta,n)$.
\end{theorem}

\proof Under the conditions of
Theorem~\ref{thm:pt}, using Theorem~\ref{erdsosind} and
Construction~\ref{LowerRT}, we have the following inequality:
$$ \RTS\left(K_\pa, \oq{\Ta}\right)\le \    \frac{1}{2}\left(1-\frac{t}{s-1}\right)<
\frac{1}{2}\left(1-\frac{1}{r}\right)\le \RTI\left(K_\pa, \Q(\Ta , n)\right). $$ Trivially,
$\frac{\Q(\Ta ,n)}{\w{n}}\le (1+o(1))\oq{\Ta}$.   Therefore,
 $$ \RTS\left(K_\pa, \frac{\Q(\Ta , n)}{\w{n}}\right)\le \RTS\left(K_\pa, \oq{\Ta}\right)< \RTI(K_\pa, \Q(\Ta , n)).\proofend
 $$

We have seen that $K_5$ has a strong phase transition at $c\sqrt{n\log
  n}$ for any $c>1$.
It follows from Theorem~\ref{thm:pt} that every clique
$K_\pa$ with $\pa\ge 5$ has a phase transition at $\Q(3,n) = \Theta\left(
\sqrt{n\log n}\right)$. 
On the other hand, $\pa=9$ and $\Ta=4$ do not satisfy the condition of Theorem~\ref{thm:pt}, and it follows from Theorem~\ref{thm:pq-1}~(b) with $p=4, q=3, i=3$ and $t=4$ that $\RT(K_9, \Q(4,n)) = \RT\big(K_9, \Q(4,o(n))\big)=1/4 $, i.e., $K_9$ does not have a strong phase transition at $\Q(4, n)$. Theorem~\ref{thm:pt} also implies that for any integer $L>0$, there exists an $\pa$ such that $K_\pa$ has more than $L$ strong phase transitions. For example, if $\pa = L!+1$, then $K_\pa$ has a strong phase transition at $\Q(\Ta,n)$ for every $\Ta$ between $2$ and $L+1$.

We also study weak phase transitions.

\begin{theorem}\label{thm:wpt}
If $K_\sa$ has a phase transition from $\Q(\Ta,n)$ to $\Q(\Ta+1,n)$, then $K_\sa$ has a weak phase transition at $\Q(\Ta, n)$.
\end{theorem}

\proof
Let $r =
\left\lfloor\frac{\sa-1}{t}\right\rfloor$ and $f(n) =
n2^{-\w{n}\log^{\frac{r}{r+1} }n}.$ 
To prove Theorem~\ref{thm:wpt}, we need that
\beq{wpt1} \RTS(K_\sa ,\Q(\Ta , f)) <  \RTI (K_\sa, \Q(\Ta , n)).\eeq
We know that $K_\sa$ has a phase transition from $\Q(\Ta , n)$ to $\Q(\Ta +1, n)$, i.e., 
$ \RTS(K_\sa ,\Q(\Ta +1, n )) <  \RTI (K_\sa, \Q(\Ta , n)).$
Therefore, to prove \eqref{wpt1},
it is sufficient to show that for $n\to\infty$, we have
\beq{kss} \RTP(n, K_\sa ,
\Q(\Ta , f(n))) \le \RTP(n, K_\sa , \Q(\Ta +1,n))+\on. \eeq  We may assume $\Q(\Ta +1,n) \le \Q(\Ta , f(n))$ since
otherwise we immediately have \eqref{kss}.
Then, by
$rt+1\le \sa \le t(r+1)$, we 
can use Construction~\ref{LowerRT} with $r$ and $t$
 as above and Theorem~\ref{thm:Kpq}~(a) with $p=t$ and $q=r+1$ to  obtain that
 $$\frac{r-1}{2r}\le\RTI(K_\sa, \Q(\Ta +1,n)) \le \RTS(K_\sa, \Q(\Ta ,
   f))\le \frac{1}{2}\left(1-\frac{1}{r}\right).$$ Hence $\RTS(K_\sa, \Q(\Ta , f)) =
   \RTI(K_\sa, \Q(\Ta +1,n))$, proving \eqref{kss}.
\proofend

We would like to have a similar result for strong phase
transitions. Unfortunately, we can prove it
  only by assuming some conditions on Ramsey numbers.  
Many of our results depend on Conjecture~\ref{conj:Qsn} and analogous conjectures. For example, \eqref{Q} and \eqref{eq:Qsn} imply that there exists a $\vartheta'$ such that
\[
\R(\ell - 1, \Q(\ell, n)) \le n^{1-\vartheta'}.
\]




If Conjecture~\ref{conj:Qsn} (b) is true for $\ell = \Ta$, then we can determine
$\RTP(n, K_{s}, \Q(\Ta,n))$. Our next result is an analogue of
Theorem~\ref{ES:RT}.  

\begin{theorem}\label{thm:rtq} 

If $r = \lfloor \frac{\sa-1}{\Ta-1} \rfloor$ and Conjecture~\ref{conj:Qsn} (b)
is true for $\ell = \Ta$, then
 $$\RT( K_\sa, \Q(\Ta,n)) =\td r.$$ 
\end{theorem}

\proof
Let $p =t-1$ and $q = r+1$. Note that
 $p(q-1)+1 \le \sa \le pq$, so by Theorem~\ref{thm:Kpq}~(b) we get the desired result.
\proofend

We also prove an extension of Theorem~\ref{thm:wpt}.

\begin{theorem}\label{thm:pt*}
 
If $\Ta\ge 2$, Conjecture~\ref{conj:Qsn} (b) is true for 
$\ell =\Ta$ and $\Ta+1$,
 and $K_\pa$ has a phase transition from $\Q(\Ta,n)$ to
$\Q(\Ta+1,n)$, 
then $K_\pa$ has a strong phase transition at $\Q(\Ta,n)$.
\end{theorem} 

\proof
Assume $r = \lfloor (\pa-1)/\Ta\rfloor$,  so $\pa<(r+1)\Ta+1$, and therefore we have
\beq{p1}
\frac{1}{2}  \left( 1-\frac{t}{s-1} \right)   <   \frac{1}{2}  \left(1-   \frac{\Ta}{(r+1)\Ta+1-1}\right) =   \frac{1}{2}  \left( 1-\frac{1}{r+1} \right).
 \eeq
By Theorem~\ref{thm:rtq} (here our $\Ta$ is $t-1$ in Theorem~\ref{thm:rtq}), we know 
that $$\RT(K_\pa, \Q(\Ta +1, n)) =  \frac{1}{2}  \left( 1-\frac{1}{r+1} \right).$$
Then by Theorem~\ref{thm:rtq} and the condition $\RT(K_\pa, \Q(\Ta +1, n)) < \RT(K_\pa, \Q(\Ta , n)) $, we have for some $r'\ge r$ that 
\beq{p2} \RT(K_\pa, \Q(\Ta , n)) = 
 \frac{1}{2}  \left( 1-\frac{1}{r'+1} \right)\ge
  \frac{1}{2}  \left( 1-\frac{1}{r+1} \right).
 \eeq
Now combining Theorem~\ref{erdsosind}, \eqref{p1} and \eqref{p2}, we have
$$\RTS\left(K_\pa,\Q\left(\Ta,\frac{n}{\w{n}}\right)\right) \le \frac{1}{2}  \left( 1-\frac{t}{s-1} \right)
< \frac{1}{2}  \left( 1-\frac{1}{r+1} \right) \le \RT(K_\pa, \Q(\Ta , n)). $$
By definition of $\Q(\Ta , n)$, it is easy to see that $\frac{\Q(\Ta ,n)}{\w{n}}\le \oq{\Ta}$, thus
$$\RTS\left( K_\pa,  \frac{\Q(\Ta ,n)}{\w{n}} \right)\le \RTS\left(K_\pa,\Q\left(\Ta,\frac{n}{\w{n}}\right)\right) < \RT(K_\pa, \Q(\Ta , n)). \proofend$$

If Conjecture~\ref{conj:Qsn} (c) is true, then what Theorem~\ref{thm:pt*} says is that  if there is a 
drop in the Ramsey-Tur\'an density while the independence number decreases down  $n^{ \frac{1}{\Ta}+o(1)}$ to $n^{ \frac{1}{\Ta+1}+o(1)}$, then there is a drop around $n^{ \frac{1}{\Ta}+o(1)} $.


We also characterize weak phase transitions for cliques.

\begin{theorem}\label{thm:wpt*}

If Conjecture~\ref{conj:Qsn} (b) is true for $\ell = \Ta+1$ and $K_\sa$ has a phase transition from $\Q(\Ta,n)$ to $\Q(\Ta+1,n)$,  then there
exists an $\eps >0$ such that
for every $\w{n}\to\infty$ slowly,
if $\ph_{\eps}(n) = 2^{-\w{n}\log^{1-\eps} n}$,
 then $K_\sa$ has a $\ph_{\eps}$-phase-transition, i.e., weak phase transition at $\Q(\Ta,n)$, and $K_\sa$ does not have a phase transition from 
$\ph_{\eps}(n) \Q(\Ta,n)$ to $\Q(\Ta+1, n)$.
\end{theorem} 

\proof
If Conjecture~\ref{conj:Qsn} (b) is true for $\ell = t+1$, then for every $\eps>0$, we have
\beq{eq:wpt*1}
\Q(\Ta +1,n) \le \varphi_{\eps}(n) \Q\left(t,n\right)  \le \Q\left(t,\varphi_{\eps}(n)n\right),
\eeq
 where the second inequality holds by the definition of $\Q(\Ta ,n)$.
Let $r = \left\lfloor(\sa-1)/t\right\rfloor$ and $\eps = \frac{r}{r+1}$.
Using the proof of Theorem~\ref{thm:wpt} (or { Theorem~\ref{thm:Kpq}~(b)} with $p=t$ and $q=r+1$), we know that 
\beq{eq:wpt*2}
\RT(K_\sa, \Q(\Ta +1,n)) = \RT\left(K_\sa, \Q\left(t,\varphi_{\eps}(n)n\right) \right).
\eeq
Now combining~\eqref{eq:wpt*1} and \eqref{eq:wpt*2}, we have  
$$ \RT(K_\sa, \Q(\Ta +1,n)) = \RT\left(K_\sa, \varphi_{\eps}(n)\Q\left(t, n\right) \right), $$
which implies the desired result.
\proofend


%
%
%
%
%

If Conjecture~\ref{conj:Qsn} (b) is true, then the assumptions of all
Theorems and Corollaries in this section
also hold.  Under this assumption, we list $\RT\left( K_{13}, f
\right)$ in Table~\ref{pttK13}, which makes our results easier to
understand.  We have three
types of functions $f(n)$: $\Q\left(t, n \right), \oq{\Ta}$ and $\Q(\Ta , g_q(n))$ where $g_q(n) = n2^{-\w{n}\log^{1-1/q} n}$. We also provide a larger table, Table~\ref{ptt} for cliques $K_4$ to $K_{13}$ in Appendix~\ref{sec:table}.

 \begin{table}[htdp]
\begin{center}
\begin{equation*}
{
\renewcommand{\arraystretch}{2.1}
 \begin{array}{r|c|c|c||r|c|c|c}\hline
 \rowcolor{mygray3} 
 & f(n) & \RT(K_{13}, f) & \textnormal{Result from} & & f(n) & \RT(K_{13}, f) & \textnormal{Result from}\\
\rowcolor{mygray}
1&n& 11/24 & \textnormal{Tur\'an's Theorem}& 7&\Q\left(5, \nq\right)&  1/3 &  \textnormal{ Theorem~\ref{thm:pq-1}~(b)}   \\
\rowcolor{mygray2}
2&o(n)  & 5/12  & \textnormal{Construction~\ref{LowerRT}}& 8&\oqd{5}  &\le 7/24   &  \textnormal{  Theorem~\ref{erdsosind}}  \\
\rowcolor{mygray}
3&\Q\left(3, \nq\right) & 5/12  & \textnormal{Construction~\ref{LowerRT}} & 9&\Q(5,g_2(n)) & 1/4  & \textnormal{ Theorem~\ref{thm:pq-1}~(b) }\\
\rowcolor{mygray2}
4&o\left(\sqrt{n\log n}\right)  &3/8 & \textnormal{Corollary~\ref{Cor:erdsosind} } & 10&\Q\left(7, \nq\right)&   1/4 & \textnormal{Theorem~\ref{thm:pq-1}~(b) } \\
\rowcolor{mygray}
5&\Q\left(4, \nq\right)&  3/8 & \textnormal{Theorem~\ref{thm:pq-1}~(b)} & 11&\oqd{7}& 0 &   \textnormal{Theorem~\ref{thm:pq-1}~(a)} \\
\rowcolor{mygray2} 
6&\oqd{4}& 1/3 & \textnormal{Corollary~\ref{Cor:erdsosind}}  \\
 \end{array} 
 }
 \end{equation*} \caption{Phase Transitions for $K_{13}$}
 \label{pttK13}
 \end{center}
\end{table}

\section{Tools}\label{sec:tools}

The method of Dependent Random Choice was developed by F\"uredi, Gowers,
Kostochka, R\"odl, Sudakov, and possibly many others.
The next lemma is taken from Alon, Krivelevich and Sudakov \cite{AlonKrSu}.
Interested readers may check the survey paper on this method by Fox and Sudakov~\cite{fs_survey}.

\begin{lemma}{(Dependent Random Choice Lemma)}\label{deprc}
Let $a,d,m,n,r$ be positive integers. Let $G=(V,E)$ be a graph with
$n$ vertices and average degree $d=2e(G)/n$. If there is a
positive integer $t$ such that
\begin{equation}\label{drc}
\frac{d^t}{n^{t-1}} - \binom{n}{r}\left(\frac{m}{n}\right)^t \ge a,
\end{equation}
then $G$ contains a subset $U$ of at least $a$ vertices such that every $r$ vertices in $U$ have at least $m$ common neighbors.
\end{lemma}

Conlon, Fox, and Sudakov~\cite{CFS} extended Lemma \ref{deprc} to hypergraphs. The {\em weight} $w(S)$ of a set $S$ of edges in a hypergraph is the number of vertices in the union of these edges. 

\begin{lemma}{(Hypergraph Dependent Random Choice Lemma).}\label{hydeprc}

Suppose $s, \Delta$ are positive integers, $\eps, \delta > 0$, and
$G_r = (V_1, \ldots, V_r; E)$ is an $r$-uniform $r$-partite hypergraph
with $|V_1|=\ldots=|V_r| = N$ and at least $\eps N^r$ edges. Then
there exists an $(r-1)$-uniform $(r-1)$-partite hypergraph $G_{r-1}$
on the vertex sets $V_2,\ldots,V_r$ which has at least
$\frac{\eps^s}{2}N^{r-1}$ edges and such that for each nonnegative
integer $w\le (r-1)\Delta$, there are at most $4r\Delta
\eps^{-s}\beta^s w^{r\Delta}r^wN^w$ dangerous sets of edges of
$G_{r-1}$ with weight $w$, where a set $S$ of edges of $G_{r-1}$ is
dangerous if $|S|\le \Delta$ and the number of vertices $v\in V_1$
such that for every edge $e\in S, e+v\in G_r$ is less than $\beta N$.
\end{lemma}


\section{Proofs of Theorems~\ref{thm:ramsey_sparse} and~\ref{sharp_k5}}\label{sec:proof_drc}

%

In this section we use Lemma~\ref{deprc} to prove Theorems~\ref{thm:ramsey_sparse} and~\ref{sharp_k5}.

\Proof of Theorem~\ref{thm:ramsey_sparse}.
The general bound \eqref{eq:ramsey} on Ramsey numbers implies that there exists a constant $\vartheta>0$ (depending on $\ell$ and $c$) such that $\R\left(s, \Q(\ell,n)n^c \right) <  n^{1-\vartheta}$. Assume that $G=G_{n,\ell}$ has more than $\eps n^2$ edges and $\eps > n^{-\vartheta^2/2s}$. We apply Lemma~\ref{deprc} to $G$ with
\begin{align*}
&r= s,\qquad d=2\eps n,\qquad t=2s/\vartheta\text{and}a = m =  \R(s, \Q(\ell,n)n^c).
\end{align*}

Now the condition of Lemma~\ref{deprc}, \eqref{drc} is satisfied as 
$$\frac{d^t}{n^{t-1}}-\binom{n}{r}\frac{m^t}{n^t} >
(2\eps)^t n-n^s\cdot n^{-\vartheta\cdot 2s/\vartheta}>\eps^t n > n^{1-\frac{\vartheta^2}{2s}\cdot\frac{2s}{\vartheta}}>a.$$ Therefore we can
use Lemma~\ref{deprc} (with the parameters $a, d, m, r, t$ as
above) to find a  subset $U$ of the vertices of  $G$ with $|U| = a$ such that all
subsets of $U$ of size $r$ have at least $m$ common neighbors.  The
set $U$ does not contain an independent set of size $\Q(\ell,n)n^c$, so $H_{n,\ell}[U]$ contains a $K_s$. Denote by $W$ the common
neighborhood of the vertices of this $K_s$. It follows that $|W|\ge m$. Then  $H_{n,\ell}[W]$  also contains a $K_{s}$, which together with the
$K_s$ found in $H_{n,\ell}[U]$ forms a $K_{2s}$.
\proofend

\Proof of Theorem~\ref{sharp_k5}.
Let $\eps_n = \w{n}^{-1/4}$.
Assume that there is a $K_5$-free graph $\Gn$ with 
\beq{AssuA}e(\Gn)\ge\eps_n n^2\text{and}\alpha(\Gn) < \frac{\sqrt{n\log n}}{\w{n}}. \eeq
We apply Lemma~\ref{deprc} to $\Gn$ with
$$a=\frac{4n}{\w{n}^2},\qquad  r=3,\qquad d=2\eps_n n,\qquad m=\sqrt{n\log n}\qquad\textnormal{and}\qquad t=7.$$
Now the condition of Lemma~\ref{deprc}, \eqref{drc} is satisfied as
$$ \frac{d^t}{n^{t-1}} - \binom{n}{r}\left(\frac{m}{n}\right)^t \ge (2\eps_n)^{7} n -n^3 \left(\frac {\log n}{n}\right)^{7/2}>  \eps_n^7 n \ge \frac{n}{\w{n}^{7/4}}>a. $$

So there exists a vertex subset $U$ of $G$ with $|U| = a=4n/\w{n}^2$
such that all subsets of $U$ of size $3$ have at least $m$ common
neighbors.  It follows from~\eqref{beta} that either $U$ has an independent set of size at least
$\left(\frac{1}{ \sqrt{2} }-o(1)\right)\sqrt{\frac{4n}{\w{n}^2}
  \log\left(\frac{4n}{\w{n}^2}\right)}> \alpha(\Gn)$, or $\Gn[U]$ contains a
triangle. 
In the latter case, 
denote by $W$ the common neighborhood of
the vertices of this triangle. It follows that $|W|\ge m=\sqrt{n\log n} > \alpha(\Gn) $,  so $\Gn[W]$ contains an edge, and this edge forms a $K_5$ with the triangle.
\proofend

\remark It is interesting to note that similar proof could be obtained with the following values:
$$\eps_n = \frac{\log\log\big(\w{n}/2\big)}{\log\big(\w{n}/2\big)},~ a=\frac{\sqrt{n\log n}}{\w{n}},~  r=2,~ d=2\eps_n n,~ m=\frac{4n}{ \w{n}^2}~\textnormal{and}~ t=\frac{\log  n}{\log\big(\w{n}/2\big)}.$$ In the proof of Theorem~\ref{sharp_k5}, we find a triangle in $U$ and then find an edge in the common neighborhood of vertices of that triangle. For those  new values, we find an edge in $U$ and then find a triangle in the common neighborhood of vertices of that edge.

\section{Proofs of Theorems~\ref{thm:Kpq}~and~\ref{thm:K2t}}\label{sec:proof_hydrc}

The proofs of Theorems \ref{thm:pq-1} and \ref{thm:Kpq} 
are very similar, therefore the proof of  Theorem \ref{thm:pq-1} is put into Appendix~\ref{sec:A1}.


We start by sketching the proof of Theorem~\ref{thm:Kpq} (a).  Suppose that $\Gn$
has more than $\left(\frac{q-2}{q-1}+\delta\right)\frac{n^2}{2}$ edges
and is $K_{pq}$-free, then we apply Szemer\'edi's Regularity Lemma to $\Gn$ and
find a $K_q$ in the cluster graph $\Rk$  (see below). 
Let $V_1,\ldots, V_q$ be the vertices of a $q$-clique in the cluster graph. We use Lemma~\ref{deprc} to find a $K_{2p}$ in $V_{q-1}\cup V_{q}$ and use Lemma~\ref{hydeprc} to find a $K_p$ in each $V_i$ for $1\le i \le q-2$ such that 
these cliques together form a $K_{pq}$ in $\Gn$. 
 The details are below.

 
\Proof of Theorem~\ref{thm:Kpq}. First we prove Theorem~\ref{thm:Kpq}~(a).
  Suppose to the contrary that there is a $K_{pq}$-free graph $\Gn$ with $n$ sufficiently large, 
  $$e(\Gn) \ge \left(\frac{q-2}{q-1}+\delta\right)\frac{n^2}{2}\text{and}\alpha(\Gn)<\Q\left(p,f(n)\right).$$

We apply Szemer\'edi's Regularity Lemma to $\Gn$ with regularity
parameter $\rho = \delta/2^{2^q}$ to get a cluster graph $\Rk$ on $k$
vertices where the vertices of $\Rk$ are the clusters of the
  Szemer\'edi Partition, and  
adjacent if the pair is $\rho$-regular and has density at least
$\delta/2$.  It is standard to check that the number of edges of $\Rk$
is at least $\left( \frac{q-2}{q-1} +
\frac{\delta}{2}\right)\frac{k^2}{2}$. So, by Tur\'an's Theorem $\Rk$
contains a $K_q$, and by Claim~\ref{clusterQ}, we can find a $K_{pq}$
in $\Gn$, a contradiction.

 To complete the proof, it 
 is sufficient to prove 
the following assertion. 

 \begin{claim}  \label{clusterQ}
 If $\alpha(\Gn)<\Q\left(p, n2^{-\w{n}\log^{1-1/q}n}\right)$
and there exists a $K_q$ in a cluster graph of $\Gn$, then we can find
a $K_{pq}$ in $\Gn$. 
 \end{claim}

 There exist $q$ vertices in $\Rk$, denoted by $V_1,\ldots,V_q$, that induce a $K_q$.  We
 define a $q$-uniform $q$-partite hypergraph $H^0$ whose vertex set is
 $\bigcup V_i$ and edge set $E(H^0)$ is the family of $q$-sets that
 span $q$-cliques in $\Gn$ and contain one vertex from each of
 $V_1,\ldots,V_q$. Let $N = |V_i| = n/k$, then by the counting lemma,
 $|E(H^0)|\ge \eps_0 N^q$, where $\eps_0 > \left( \delta/3
 \right)^{\binom{q}{2}}$.
Let $$\beta = f(n)/N,\quad
s = \log^{\frac{1}{q} }n, \quad \eps_i = \eps_0^{\log^{\frac{i}{q} }n}2^{-\frac{s^i-1}{s-1}}, \quad r_i = q-i, \quad \Delta_i = p^{r_i} \quad \textnormal{and}\quad w_i = pr_i.$$ 
We start from $H^0$. For $1\le i \le q-2$ we  apply Lemma \ref{hydeprc} to $H^{i-1}$ with $\Delta = \Delta_i, \eps = \eps_{i-1}, r = r_{i-1}$ and $w = w_i$ to get $H^i$.
Note that $\Delta, \eps_0, r, w$ and $k$ are all constants. 
It is easy to check that for $1\le i\le q-2$, we have
\begin{eqnarray*}
4r\Delta\eps^{-s}\beta^s w^{r\Delta}r^{w}N^{w}
 &=&O\left(2^{2\log^{\frac{i-1}{q}}n}\eps_0^{-\log^{\frac{i}{q} }n}k^{  \log^{\frac{1}{q} }n }2^{-\w{n}\log n} N^{w}\right)\\
 &=&O(n^{-\w{n}/2})
 = o(1) < 1.
\end{eqnarray*}
Then by Lemma \ref{hydeprc} there exists an $r_i$-uniform
$r_i$-partite hypergraph $H^i$ on the vertex sets $V_{i+1},\ldots,V_q$
that contains at least $\eps_i N^{r_i}$ edges and contains no
dangerous sets of $\Delta_i$ edges on $w_i$ vertices. (Recall that a set $S$ of $\Delta_i$ edges on $w_i$ vertices is dangerous if the number of vertices $v \in V_i$ for which for every edge $e\in S, e+v \in H^{i-1}$ is less than $\beta N$). Now we have a
hypergraph sequence $\{H^\ell\}_{\ell=0}^{q-2}$.  We will prove by
induction on $i$ that there is a $p$-set $A^{q-\ell}\subset V_{q-\ell}$ for
$0\le \ell\le i$ such that $\Gn\left[A^{q-\ell}\right]=K_p$ and
$H^{q-i-1}\left[\bigcup_{\ell=0}^{i} A^{q-\ell}\right]$ is complete $r_{q-i-1}$-partite. Note that if a vertex set $T$ is an edge
of $H^0$, then $\Gn[T]$ is a $q$-clique. So
$\Gn\left[\bigcup_{\ell=0}^{q-1} A^{q-\ell}\right] = K_{pq}$, which
will prove Claim~\ref{clusterQ}.

We first show that the induction hypothesis holds for $i=1$. Note that $r_{q-2}$ = 2, so $H^{q-2}$ is a bipartite graph on $2N$ vertices with at least $\eps_{q-2}N^2$ edges. We now apply Lemma~\ref{deprc} to $H^{q-2}$ with 
$$a=2\beta N,  \qquad d=\eps_{q-2}N, \qquad t=s, 
  \qquad r=p \qquad \textnormal{and} \qquad m=\beta N.$$ 
We check condition \eqref{drc}:
\begin{eqnarray*}
\frac{(\eps_{q-2}N)^s}{(2N)^{s-1}}-{\binom{2N}{p}}\left( \frac{\beta N}{2N} \right)^s 
&\ge& (\eps_0/2)^{\log^{1-1/q}n}N - n^p k^s2^{-\w{n}\log n} \\
&=&  (\eps_0/2)^{\log^{1-1/q}n}N - o(1)
\ge 2\beta N.
\end{eqnarray*}
Therefore we have a subset $U$ of $V_{q-1}\cup V_{q}$ with $|U| = 2\beta N$ such that 
every $p$ vertices in $U$ have at least $\beta N$ common neighbors in $H^{q-2}$. Either $V_{q-1}$ or $V_{q}$ contains at least half of the vertices of $U$, so w.l.o.g. we may assume 
that $U' = U\cap V_{q-1}$ contains at least $\beta N = m$ vertices. 
Because $\alpha(\Gn) < \Q(p,m)$, the vertex set $U'$ contains a $p$-vertex set $A^{q-1}$
such that $\Gn\left[A^{q-1}\right] = K_p$. The vertices of $A^{q-1}$ have at least $m$ common neighbors in $V_q$, so their common neighborhood also contains a $p$-vertex subset $A^q$ of $V_q$ such that $\Gn[A^q] = K_p$.
 Now $H^{q-2}\left[A^{q-1}\cup A^q\right]$ is complete bipartite.
 We are done with the base case $i=1$.

 For the induction step, assume that the induction hypothesis
holds for $i-1$, then we can find a complete $r_{q-i}$-partite subhypergraph  $\widetilde{H}^{q-i}$  of $H^{q-i}$ spanned by $\bigcup_{\ell=0}^{i-1} A^{q-\ell}$, where $\Gn[A^{q-\ell}] = K_p$ for  every $\ell$. The hypergraph
$H^{q-i}$ has no dangerous set of $\Delta_{q-i}$ edges on $w_{q-i}$ vertices, and $\widetilde{H}^{q-i}$ contains $pi = w_{q-i}$ vertices and $p^i = \Delta_{q-i}$ edges,
so $\widetilde{H}^{q-i}$  is not dangerous.
 Then we can find a set $B$ of $\beta N$ vertices in $V_{q-i}$ such that 
 for every edge $e\in \widetilde{H}^{q-i}$ and every vertex $v\in B$, $e+v \in H^{q-i-1}$, which means that $ H^{q-i-1}\left[B\cup\bigcup_{\ell=0}^{i-1} A^{q-\ell}\right]$ is complete $r_{q-i-1}$-partite. Then,
 because $\alpha(\Gn) < \Q(p,\beta N)$, we can find a $p$-vertex subset $A^{q-i}$ of $B$ such that $\Gn[A^{q-i}] = K_p$.

We apply Theorem~\ref{thm:Kpq}~(a) to get the upper bound in Theorem~\ref{thm:Kpq}~(b). Let $r=q-1$, which implies
 $rt+1\le pq-i$, so the lower bound is realized by \eqref{lowerQ}  with the parameters $r,t$ as above: 
\[
 \RTI\left(K_{pq-i}, \Q\left(t+1, \frac{n}{r}\right)\right) \ge  \RTI\left(K_{rt+1}, \Q\left(t+1, \frac{n}{r}\right)\right) 
\ge \frac{1}{2}\left(1-\frac{1}{q-1}\right). \proofend
\]


The proof of Theorem~\ref{thm:K2t} is a combination of Claim~\ref{clusterQ} and an easy application of Szemer\'edi's Regularity Lemma,
(see the Appendix of Balogh-Lenz~\cite{BL2} 
for similar proofs). The idea is  that instead of proving only that the cluster graph is $K_q$-free, like in the proof of Theorem~\ref{thm:Kpq}, we also bound the density of regular pairs.

\Proof of Theorem~\ref{thm:K2t}.
Given $\eps>0$, let $\rho = \eps/2^{2^t}$ and $M=M(\rho)>1/\rho$ be the upper bound on the number of partitions guaranteed by Szemer\'edi's Regularity Lemma
with regularity parameter $\rho$.
  Suppose we have a $K_{2t}$-free graph $\Gn$ with 
  $$e(\Gn)\ge\left(\frac{(t-1)(q-2)}{t(q-1)}+\eps\right)\frac{n^2}{2}\text{and}\alpha(\Gn) < \Q\left(t,\frac{\eps n}{M}\right).$$ 
  We apply Szemer\'edi's Regularity Lemma to $\Gn$ with regularity
  parameter $\rho$ to get a cluster graph $\Rk$ on $k\le M$ vertices
  where two vertices are adjacent if the pair is $\rho$-regular and
  has density at least $\eps/2$.  It is standard to check that more
  than $\left( \frac{(t-1)(q-2)}{t(q-1)} +
  \frac{\eps}{2}\right)\frac{n^2}{2}$ edges of $\Gn$ are between pairs
  of classes that are $\rho$-regular and have density at least
  $\eps/2$.

Assume that the density $d$ of a $\rho$-regular 
 pair $(V_i, V_j)$ is at least $\frac{t-1}{t}+\eps$. Because $\alpha(\Gn) < \Q(\Ta ,\eps  n/M)$ and $|V_i| \ge \eps  n/M$, there is  a $t$-clique in $V_i$, each of whose vertices has at least $(d-\rho)|V_j|\ge (\frac{t-1}{t}+\frac{\eps}{2})|V_j|$ neighbors in $V_j$,
 hence vertices of this $t$-clique have 
at least $\eps |V_j| $ common neighbors. Then we can find a $t$-clique in their common neighborhood since $\alpha(\Gn) < \Q(\Ta ,\eps  n/M)$ and $\eps |V_j| \ge \eps  n/M$. Thus we find a $K_{2t}$ in $\Gn$, a contradiction. Therefore the density of any $\rho$-regular pair is at most  $\frac{t-1}{t}+\eps$. Then $\Rk$ has at least  
$$\left.\left(\frac{(t-1)(q-2)}{t(q-1)}+\frac{\eps}{2}\right)\frac{n^2}{2} \cdot \left(\left(\frac{t-1}{t}+\eps\right)\left(\frac{n}{k}\right)^2\right)^{-1}\right. >  \left(\frac{q-2}{q-1}+\frac{\eps}{4}\right)\frac{k^2}{2}$$ 
edges, so there is  a $K_q$ in $\Rk$. Then, by Claim~\ref{clusterQ}, there is  a $K_{pq}$ in $\Gn$,  a contradiction.
\proofend


%


{

\section{Open problems}\label{sec:last}

We proved that $\RT\left(K_5,o\left(\sqrt {n\log n}\right)\right)=0$, and
it was known that $\RT\left( K_5 , \Q(3, n/2)
\right) = 1/4 $. It would be interesting to know if
there is any sharper transition  at  $c\cdot \Q(3, n/2)$ for $c<1$,
hence it is natural to propose the following two problems:

\Question{
Determine
$\RTP\left(n, K_5,(1-\eps)\Q(3, n/2)\right)$.}

\Question{
Determine $\RTP\left(n, K_5, c\cdot \Q(3, n/2)\right)$  for $0<c<1$.}

We proved that if Conjecture~\ref{conj:Qsn} (b) is true, then $\RTS\big(
K_{2\Ta}, o\big(\Q(\Ta ,n)\big) \big) \le \frac{\Ta-1}{4\Ta}. $ Note that the
Bollob\'as-Erd\H os graph gave matching lower bound for $\Ta=2$, so
finding constructions to give matching lower bounds on $\RTI\big(
K_{2\Ta}, o\big(\Q(\Ta ,n)\big)\big) $ is a very challenging problem.
The most interesting case is $K_6$. A construction improving the lower bound on $\RTI\big(K_6,$ $o\left(\sqrt{n\log n}\right)\big)$ would imply several improvements in the spirit
 of Construction~\ref{LowerBE}, i.e., 
 then we could
 replace the Bollob\'as-Erd\H os graph in Construction~\ref{LowerBE}  
 with this new construction for $K_6$, replace $H_{k,3}$ with $H_{k,4}$, 
 and then optimize the class sizes.
 Probably such a  construction, if exists,  is an
extension of the Bollob\'as-Erd\H os graph.  There are generalizations
of the Bollob\'as-Erd\H os graph in~\cite{BL1, BL2, EHSSSzA, EHSSSzB}.
The case of $K_6$  is asked below:
  
  \Question{
  Determine $\RTS\left(K_6,o\left(\sqrt{n\log n}\right)\right)$ and $\RTI\left(K_6,o\left(\sqrt{n\log n}\right)\right)$.}
  
We have $1/6$ as an upper bound.  Sudakov proved that
$\RT\left( K_6, f \right)= 0$ for $f(n) = \Q\left(3, 
n2^{-\w{n}\sqrt{\log n} }\right)$, but it is not clear what happens
when $f(n)$ is between $\Q\left(3, n2^{-\w{n}\sqrt{\log n}} \right)$ and
$o\left(\sqrt{n\log n}\right)$. In particular, we would like to know
the answer to the following question:

} 

} 

    \Question{ For  which function $f(n)$ does $K_6$ have a strong phase transition to $0$, i.e.,
      $0 = \RTS\left(K_6, o(f)\right)< \RTI\left(K_6, f\right)$?}
      
One surprising phenomenon is that $\RT(K_4,o(\sqrt {n\log
  n}))=0=\RT(K_5, o(\sqrt { n\log n } ) ).$
  We know that $\RTS(K_6,o(\sqrt {n\log n}))\le 1/6 < 1/4 = \RT(K_7, o(\sqrt {n\log n})).$
   It would be interesting to know
if $$\RT\left(K_7,o\left(\sqrt {n\log n}\right)\right) = \RTS\left(K_8,o\left(\sqrt
{n\log n}\right)\right).$$

\appendices
\section{Proof of Theorem~\ref{thm:pq-1}}\label{sec:A1}
This proof is very similar to 
that
of Theorem~\ref{thm:Kpq} in Section~\ref{sec:proof_hydrc},
so  we skip some details.
We first prove Theorem~\ref{thm:pq-1}~(a).
  Suppose to the contrary that there is a $K_{pq-1}$-free graph $\Gn$ with 
  $n$ sufficiently large, 
  $$e(\Gn)\ge\left(1- \frac{1}{q-1}+\delta\right)\frac{n^2}
  {2}\text{and}\alpha(\Gn) < \Q\left(p, f(n) \right).$$ 
  
  Just as what we did in the proof of Theorem~\ref{thm:Kpq}, we apply Szemer\'edi's
  Regularity Lemma to $\Gn$ with regularity parameter $\rho = \delta/2^{2^q}$ to get
  a cluster graph $R$ on $k$ vertices.  Similarly to the proof of
  Theorem~\ref{thm:Kpq}, 
  we can find $q$ vertices $V_1,\ldots,V_q$ that span a $K_q$ in $R$.
Now consider a $q$-uniform $q$-partite hypergraph $H^0$ whose vertex set is $\bigcup V_i$ and edge set $E(H^0)$ is the 
family of $q$-sets that span $q$-cliques in $\Gn$ and contain one vertex from each of $V_1,\ldots,V_q$. Let $N = |V_i|  = n/k$, then  by the counting lemma, $|E(H^0)|\ge \eps_0 N^q$, where $\eps_0 > \left( \delta/3 \right)^{\binom{q}{2}}$.
Let 
\[
\beta =f(n)/N, ~ s = \log^{\frac{1}{q-1}} n, ~ \eps_i = \eps_0^{s^i}/2^{\frac{s^i-1}{s-1}},~
r_i = q-i, ~ w_i = pr_i-1~\textnormal{and}~\Delta_i = p^{r_i-1}(p-1).
\]
We start from $H^0$. For $1\le i\le q-2$ we apply Lemma \ref{hydeprc} to $H^{i-1}$ with $\Delta = \Delta_i, \eps = \eps_{i-1}, r = r_{i-1}$ and $w = w_i$ to get $H^i$.
It is easy to check that for $1\le i\le q-2$, we have 
\begin{eqnarray*}
4r\Delta \eps^{-s}\beta^s w^{r\Delta}r^{w}N^{w}
 &=&O\left(2^{2\log^{\frac{i-1}{q-1} }n}\eps_0^{-\log^{\frac{i}{q-1} }n}k^s2^{-\w{n}\log n} N^{w}\right)\\
 &=&O(n^{-\w{n}/2})
 = o(1) < 1.
\end{eqnarray*}
Then, by Lemma \ref{hydeprc}, there exists an $r_i$-uniform $r_i$-partite hypergraph $H^i$
on the vertex sets $V_{i+1},\ldots,V_q$ that contains
 at least $\eps_i N^{r_i}$ edges and contains no dangerous set of $\Delta_i$ edges on $w_i$ vertices.

Note that $r_{q-2}$ = 2, so $H^{q-2}$ is a bipartite graph on $2N$ vertices with at least $\eps_{q-2}N^2$ edges. We now apply Lemma \ref{deprc} to $H^{q-2}$ with 
\begin{align*}
a=2\beta N,\qquad d=\eps_{q-2}N, \qquad t=2p/\vartheta , \qquad 
 r=p \qquad
\textnormal{and}\qquad m=\R\left(p-1, \Q\left(p, f(n)\right)\right).
\end{align*}
Note that $m< n^{1-\vartheta}$. We check condition \eqref{drc}:
\begin{eqnarray*}
\frac{(\eps_{q-2}N)^t}{(2N)^{t-1}}-{\binom{2N}{p}}\left( \frac{m}{2N} \right)^t 
&\ge& (\eps_0/2)^{t\log^{\frac{q-2}{q-1}}n}N - n^p k^t n^{-2p} \\
&=&  (\eps_0/2)^{t\log^{\frac{q-2}{q-1}}n}N - o(1) \ge a.
\end{eqnarray*}
Therefore we have a subset $U$ of $V_{q-1}\cup V_{q}$ with $|U| =
2\beta N$ such that every $p$ vertices in $U$ have at least
$m$ common neighbors in $H^{q-2}$. Either $V_{q-1}$ or $V_{q}$
contains at least half of the
vertices of $U$, so w.l.o.g. we may assume that
$U'
:= U\cap V_{q-1}$ contains at least $\beta N$ vertices. Since
$\alpha(\Gn) < \Q(p,\beta N)$, the vertex set $U'$ contains a $p$-vertex set
$A^{q-1}$ such that $\Gn\left[A^{q-1}\right] = K_p$. The vertices of $A^{q-1}$ have at least 
$m=\R\left(p-1, \Q\left(p, \beta N\right)\right)$ common neighbors in $V_q$, so their common neighborhood 
contains a $(p-1)$-vertex subset $A^q$ of $V_q$ such that $\Gn[A^q] = K_{p-1}$. Now
$H^{q-2}\left[A^{q-1}\cup A^q\right]$ is complete bipartite.
Then, similarly to the proof of Theorem~\ref{thm:Kpq}, for $1\le i\le
q$, we can find a subset $A^{i}$ of $V_{i}$ satisfying the following conditions:
\begin{itemize}
\item $\Gn[A^q] = K_{p-1}$.
\item For $1\le i< q$, $\Gn[A^{i}]=K_p$.
\item $H^0[\bigcup_{i=1}^{q} A^{i}]$ is complete $q$-partite.
\end{itemize}
 If a vertex set $T$ is an edge of $H^0$, then $\Gn[T] = K_q$. So $\Gn\left[\bigcup_{i=1}^{q} A^{i}\right] = K_{pq-1}$, which is a contradiction.

For Theorem~\ref{thm:pq-1}~(b), the upper bound is obvious from Theorem~\ref{thm:pq-1}~(a). Let $r=q-1$, then $rt+1\le pq-i$, so the lower bound is realized by \eqref{lowerQ} with the parameters $r,t$ as above:
$$ \RTI\left(K_{pq-i}, \Q\left(t+1, \frac{n}{q-1}\right) \right)
\ge \RTI\left(K_{rt+1}, \Q\left(t+1, \frac{n}{r}\right) \right)\ge\frac{1}{2}\left(1- \frac{1}{q-1}\right). \proofend
$$



\section{}\label{sec:table}
Assuming Conjecture~\ref{conj:Qsn} is true,
we summarize our
results in Section~\ref{sec:results} by listing $\RT\left( K_\sa, f
\right)$
for $\sa\le 13$ in Table~\ref{ptt}. 

Note that under the assumption that Conjecture~\ref{conj:Qsn} (b) is true,
our results for $f(n) = \oq{t}$
can be viewed as results on $\RT\left( K_\pa, o\left(\Q(\Ta ,n)\right) \right)$.
Conjecture~\ref{conj:Qsn} is true for 
$\ell=3,4$,
 therefore, our results on $K_4,\ldots,K_8$ and results in Row~$1$ to Row~$7$ do not depend on Conjecture~\ref{conj:Qsn}.
 
 
\begin{table}[htdp]
\begin{center}
\begin{equation*}
{
\renewcommand{\arraystretch}{2.1}
\begin{array}{r|c|c|c|c|c|c||c|c|c|c|c}\hline 
\rowcolor{mygray3}
&&K_4 & K_5 & K_6 & K_7 & K_8 & K_9 & K_{10}& K_{11} &  K_{12} &K_{13}  \\\hline 
\rowcolor{mygray}
1&n& \dfrac{1}{3} &  \dfrac{3}{8} & \dfrac{2}{5} & \dfrac{5}{12} & \dfrac{3}{7} & \dfrac{7}{16} & \dfrac{4}{9} & \dfrac{9}{20} & \dfrac{5}{11}  & \dfrac{11}{24} \\
\rowcolor{mygray2}
2&o(n) &\dfrac{1}{8} &\dfrac{1}{4} & \dfrac{2}{7} & \dfrac{1}{3} & \dfrac{7}{20} & \dfrac{3}{8} & \dfrac{5}{13} & \dfrac{2}{5} & \dfrac{13}{32}  & \dfrac{ 5}{12}   \\
\rowcolor{mygray3}
3&g_q(n) &2: 0 &  & 3: \dfrac{1}{4} &  & 4: \dfrac{1}{3} &  & 5: \dfrac{3}{8} &  & 6: \dfrac{2}{5}  &  \\\hline
\rowcolor{mygray}
4&\Q\left(3, \nq\right)&	 &\dfrac{1}{4}  & \dfrac{1}{4} & \dfrac{1}{3} &  \dfrac{1}{3} & \dfrac{3}{8} & \dfrac{3}{8} &  \dfrac{2}{5} &  \dfrac{2}{5} & \dfrac{5}{12}  \\
\rowcolor{mygray2}
5&o\left(\sqrt{n\log n}\right)& & 0 & \le \dfrac{1}{6} & \dfrac{1}{4} & \le \dfrac{2}{7} & \le \dfrac{5}{16} & \dfrac{1}{3} &\le \dfrac{7}{20 } &  \le \dfrac{8}{22 }   &\dfrac{3}{8} \\
\rowcolor{mygray3}
6&\Q(3, g_q(n))& &  & 2: 0 &  & 2: \dfrac{1}{4} & 3: \dfrac{1}{4} &  & 3: \dfrac{1}{3} &  4:  \dfrac{1}{3}& \\\hline 
\rowcolor{mygray}
7&\Q\left(4, \nq\right)& &  &  & \dfrac{1}{4} & \dfrac{1}{4} & & \dfrac{ 1}{3} & \dfrac{1}{3} & \dfrac{1}{3}  &  \dfrac{3}{8} \\\hline
\hline
\rowcolor{mygray2} 
8&\oqd{4}& &  &  & 0 & \le \dfrac{3}{16} &  & \le \dfrac{5}{18} & \le \dfrac{3}{10} & \le \dfrac{7}{22}   & \dfrac{1}{3}  \\
\rowcolor{mygray3}
9&\Q(4,g_q(n))& &  &  &  & 2: 0 &  & 2: \dfrac{1}{4 }& 2: \dfrac{1}{4 }&  3: \dfrac{1}{4 } &  \\\hline
\rowcolor{mygray} 
10&\Q\left(5, \nq\right)&   &  &  &  &  & \dfrac{1}{4} & \dfrac{ 1}{4} &  &  &  \dfrac{1}{3}   \\
\rowcolor{mygray2}
11&\oqd{5}&   &  &  &  &  & 0 & \le \dfrac{1}{5} & &   &\le \dfrac{7}{24}      \\
\rowcolor{mygray3}
12&\Q(5,g_q(n))&   &  &  &  &  &  & 2: 0 & &   &2: \dfrac{1}{4}    \\\hline 
\rowcolor{mygray}
13&\Q\left(6, \nq\right)&   &  &  &  &  &  & &   \dfrac{1}{4} &   \dfrac{1}{4} &  \\
\rowcolor{mygray2} 
14&\oqd{6}&   &  &  &  &  &  & & 0  & \le \dfrac{5}{24} &  \  \\
\rowcolor{mygray3}
15&\Q(6,g_q(n))&   &  &  &  &  &  & &  & 2:0  & \\\hline 
\rowcolor{mygray}
16&\Q\left(7, \nq\right)&   &  &  &  &  &  & & &  &   \dfrac{1}{4}   \\
\rowcolor{mygray2}
17&\oqd{7}&    &  &  &  &  &  & &  &  & 0    \\\hline 

\end{array}
}
\end{equation*} \caption{Phase Transitions.}
\label{ptt}
\end{center}
\end{table}
 

An entry ``$\lambda$" in the row $f(n)$ and the column $K_\sa$ means
$\RT\left( K_\sa, f \right)=\lambda$, and ``$\le \lambda$'' means
$\RTS\left( K_\sa, f\right)\le \lambda$.  In the row $\Q(\Ta , g_q(n))$,
the entries are ``$q_0:\lambda$''  meaning that
$\RT\left( K_\sa, \Q(\Ta , g_{q_0})\right)=\lambda$.

{\bf Acknowledgement:}
 We would  like to thank the referees for  careful reading of the manuscript.





\end{document}